\documentclass{article}

\PassOptionsToPackage{square, numbers}{natbib}

 \usepackage[final]{nips_2016}

 \usepackage{threeparttable}

\usepackage[utf8]{inputenc} 
\usepackage[T1]{fontenc}    
\usepackage{hyperref}       
\usepackage{url}            
\usepackage{booktabs}       
\usepackage{amsfonts}       
\usepackage{nicefrac}       
\usepackage{microtype}      

\usepackage{color}
\usepackage{color,graphicx,amssymb,amsmath,amsthm}

\usepackage{lipsum}

\newcommand\blfootnote[1]{%
  \begingroup
  \renewcommand\thefootnote{}\footnote{#1}%
  \addtocounter{footnote}{-1}%
  \endgroup
}

\newcommand{\itk}{_k}
\newcommand{\kp}{_{k+1}}
\newcommand{\eqhd}{\text{hard}}

\usepackage{amsmath}
\usepackage{etoolbox}
\usepackage{bm}
\usepackage{mathtools}
\usepackage{algorithm}
\usepackage{algorithmic}
\usepackage{subfigure}
\usepackage{xspace}
\usepackage{multirow}
\usepackage[capitalise]{cleveref}
\usepackage{array}

\usepackage{footnote}
\makesavenoteenv{tabular}
\makesavenoteenv{table}

\crefname{rl}{Rule}{Rule}

\newcommand{\MyMapTemplatePrefixc}[4]{\expandafter#1\csname#3#4\endcsname{#2{#4}}} 
\forcsvlist{\MyMapTemplatePrefixc {\def} {\mathcal}{mc}} {A,B,C,D,E,F,G,H,I,J,K,L,M,N,O,P,Q,R,S,T,U,V,W,X,Y,Z}  

\newcommand{\MyMapTemplatePrefixtb}[5]{\expandafter#1\csname#4#5\endcsname{#2{#3{#5}}}} 
\forcsvlist{\MyMapTemplatePrefixtb {\def} {\tilde}{\mathbf}{t}} {A,B,C,D,E,F,G,H,I,J,K,L,M,N,O,P,Q,R,S,T,U,V,W,X,Y,Z}  
\forcsvlist{\MyMapTemplatePrefixtb {\def} {\tilde}{\mathbf}{t}} {0,1,a,b,c,d,e,f,g,h,i,j,k,l,m,n,o,p,q,r,s,t,u,v,w,x,y,z}  

\newcommand{\MyMapTemplateNoPrefix}[3]{\expandafter#1\csname#3\endcsname{#2{#3}}}
\forcsvlist{\MyMapTemplateNoPrefix {\def} {\mathbf} } {A,B,C,D,E,F,G,H,I,J,K,L,M,N,O,P,Q,R,S,T,U,V,W,X,Y,Z}  

\def\ge{\epsilon}

\def\gl{\lambda}

\newcommand{\tabincell}[2]{\begin{tabular}{@{}#1@{}}#2\end{tabular}}

\def\bbC{{\mathbb C}}

\def\bbR{{\mathbb R}}

\def\st{\mbox{subject to~~}}

\def\eqif{\mbox{if}}

\newcommand{\grad}{\nabla}

\newcommand*{\affaddr}[1]{#1} 
\newcommand*{\affmark}[1][*]{\textsuperscript{#1}}

\title{An Empirical Study of ADMM for \\ Nonconvex Problems
}

\author{%
Zheng Xu\affmark[1],\  Soham De\affmark[1],\ M\'{a}rio A. T. Figueiredo\affmark[2],\ Christoph Studer \affmark[3],\ Tom Goldstein\affmark[1]\\
\affaddr{\affmark[1]Department of Computer Science, University of Maryland, College Park, MD}\\
\affaddr{\affmark[2]Instituto de Telecomunica\c{c}\~{o}es, Instituto Superior T\'{e}cnico, Universidade de Lisboa, Portugal}\\
\affaddr{\affmark[3]Department of Electrical and Computer Engineering, Cornell University, Ithaca, NY}\\
}

\begin{document}

\maketitle

\begin{abstract}
The alternating direction method of multipliers (ADMM) is a common optimization tool for solving constrained and non-differentiable problems. We provide an empirical study of the practical performance of ADMM on several nonconvex applications, including $\ell_0$ regularized linear regression, $\ell_0$ regularized image denoising, phase retrieval, and eigenvector computation. Our experiments suggest that ADMM performs well on a broad class of non-convex problems. Moreover, recently proposed adaptive ADMM methods,  which automatically tune penalty parameters as the method runs, can improve algorithm efficiency and solution quality compared to ADMM with a non-tuned penalty.
\end{abstract}

\blfootnote{
\hspace{-7.2mm} ZX, SD, and TG were supported by US NSF grant CCF-1535902 and by US ONR grant N00014-15-1-2676.\\
CS was supported in part by Xilinx Inc., and by the US NSF under grants ECCS-1408006 and CCF-1535897.
}

\section{Introduction}
\label{sec:intro}
The alternating direction method of multipliers (ADMM) has been applied to solve a wide range of constrained convex and nonconvex optimization problems.
ADMM decomposes complex optimization problems into sequences of simpler subproblems that are often solvable in closed form.  Furthermore, these sub-problems are often amenable to large-scale distributed computing environments \cite{goldstein2016unwrapping,taylor2016training}. ADMM solves the problem
\begin{eqnarray}
\min_{u\in \bbR^n,v\in \bbR^m}  H(u) + G(v),~~~~\st~~  Au+Bv = b, \label{eq:prob}
\end{eqnarray}
where $H:\bbR^n \! \rightarrow \! \bar{\bbR}$, $G:\bbR^m \! \rightarrow \! \bar{\bbR}$, $A\in \bbR^{p \times n}$, $B\in \bbR^{p \times m}$, and $b \in \bbR^p$,  by the following steps,
\begin{align}
u\kp =& \arg\min_{u} H(u) + \langle \gl\itk, -Au \rangle + \frac{\tau\itk}{2} \| b-Au-Bv\itk\|_2^2 \label{eq:updateu}\\
v\kp =& \arg\min_{v} G(v) + \langle \gl\itk, -Bv \rangle + \frac{\tau\itk}{2} \| b-Au\kp-Bv\|_2^2 \label{eq:updatev}\\
\gl\kp =& \gl\itk +\tau\itk (b - A u\kp -B v\kp), \label{eq:updatedual}
\end{align}
where $\gl \! \in\! \bbR^p $ is a vector of dual variables (Lagrange multipliers), and $\tau\itk$ is a scalar penalty parameter.

The convergence of the algorithm can be monitored using primal and dual ``residuals,'' both of which approach zero as the iterates become more accurate, and which are  defined as
\begin{eqnarray}
r\itk = b-Au\itk-Bv\itk,~~~~~\text{and}~~~~~d\itk = \tau\itk A^{T}B(v\itk-v^{k-1}),
\end{eqnarray}
respectively \cite{boyd2011admm}.  The iteration is generally stopped when
\begin{eqnarray}
\|r\itk\|_2  \leq \ge^{tol} \max\{\|Au\itk\|_2, \|Bv\itk\|_2, \|b\|_2 \}  ~~\text{and}~~ \|d\itk\|_2 \leq \ge^{tol} \| A^{T} \gl\itk\|_2, \label{eq:stop}
\end{eqnarray}
where $\ge^{tol} > 0$ is the stopping tolerance.

ADMM was introduced by \citet{glowinski1975approximation} and \citet{gabay1976dual}, and convergence has been proved under mild conditions for convex problems \cite{gabay1983con,eckstein1992douglas,he2015non}. The practical performance of ADMM on convex problems has been extensively studied, see \cite{boyd2011admm,goldstein2014fast,xu2016adaptive} and references therein. For nonconvex problems, the convergence of ADMM under certain assumptions are studied in ~\cite{wang2014convergence,li2015global,hong2016convergence,wang2015global}. The current weakest assumptions are given in \cite{wang2015global}, which requires a number of strict conditions on the objective, including a Lipschitz differentiable objective term. In practice, ADMM has been applied on various nonconvex problems, including nonnegative matrix factorization \cite{xu2012alternating}, $\ell_p$-norm regularization ($0<p<1$)\cite{bouaziz2013sparse,chartrand2013nonconvex}, tensor factorization \cite{liavas2015parallel,xu2016tensor}, phase retrieval \cite{wen2012alternating}, manifold optimization \cite{lai2014splitting,kovnatsky2015madmm},  random fields \cite{miksik2014distributed}, and deep neural networks \cite{taylor2016training}.

The penalty parameter $\tau\itk$ is the only free choice in ADMM, and plays an important role in the practical performance of the method. Adaptive methods have been proposed to automatically tune this parameter as the algorithm runs. The residual balancing method \cite{he2000alternating} automatically increase or decrease the penalty so that the primal and dual residuals have approximately similar magnitudes. The more recent AADMM method \cite{xu2016adaptive} uses a spectral (Barzilai-Borwein) rule for tuning the penalty parameter. These methods achieve impressive practical performance for convex problems
and are guaranteed to converge under moderate conditions (such as when adaptivity is stopped after a finite number of iterations).

In this manuscript, we study the practical performance of ADMM on several nonconvex applications, including $\ell_0$ regularized linear regression, $\ell_0$ regularized image denoising, phase retrieval, and eigenvector computation. While the convergence of these applications may (not) be guaranteed by the current theory, ADMM is one of the (popular) choices to solve these nonconvex problems. The following questions are addressed using these model problems: (i) does ADMM converge in practice, (ii) does the update order of $H(u)$ and $G(v)$ matter, (iii) is the local optimal solution good, (iv) does the penalty parameter $\tau\itk$ matter, and (v) is an adaptive penalty choice effective?

\section{Nonconvex applications}

\textbf{$\ell_0$ regularized linear regression.}
Sparse linear regression can be achieved using the non-convex, $\ell_0$ regularized problem
\begin{eqnarray}
\min_x \frac{1}{2} \| Dx - c \|_2^2 + \rho \|x\|_0, \label{eq:l0r1}
\end{eqnarray}
where $D \in \bbR^{n \times m}$ is the data matrix, $c$ is a measurement vector, and $x$ is the regression coefficients. ADMM is applied to solve problem~\eqref{eq:l0r1} using the equivalent formulation
\begin{eqnarray}
\min_{u,v} \frac{1}{2} \| Du - c \|_2^2 + \rho \|v\|_0~~~~\st~~u-v=0.
\end{eqnarray}

\textbf{$\ell_0$ regularized image denoising.}
The $\ell_0$ regularizer~\cite{dong2013efficient} can be substituted for the $\ell_1$ regularizer when computing total variation for image denoising.  This results in the formulation ~\cite{chartrand2007exact}
\begin{eqnarray}
\min_{x} \frac{1}{2} \| x-c \|_2^2 + \rho \| \grad x\|_0
\end{eqnarray}
where $c$ represents a given noisy image, $\grad$ is the linear discrete gradient operator, and $\| \cdot \|_2 /  \| \cdot \|_0$ is the $\ell_2$/$ \ell_0$ norm.
We solve the equivalent problem
\begin{eqnarray}
\min_{u,v} \frac{1}{2} \| u-c \|_2^2 + \rho \| v\|_0~~~~\st~~ \grad u-v=0.
\end{eqnarray}
The resulting ADMM sub-problems can be solved in closed form using fast Fourier transforms \cite{goldstein2009split}.

\textbf{Phase retrieval.}
Ptychographic phase retrieval~\cite{yang2011iterative,wen2012alternating} solves the problem
\begin{eqnarray}
\min_{x} \frac{1}{2}||\text{abs}(Dx)-c||_2^2,
\end{eqnarray}
where $x\in \bbC^n$, $D\in \bbC^{m \times n}$, and $\text{abs}(\cdot)$ denotes the elementwise magnitude of a complex vector. ADMM is applied to the equivalent problem
\begin{eqnarray}
\min_{u,v} \frac{1}{2}||\text{abs}(u)-c||_2^2 ~~\st u-Dv = 0.
\end{eqnarray}

\textbf{Eigenvector problem.}
The eigenvector problem is a fundamental problem in numerical linear algebra.  The leading eigenvalue of a matrix $D$ is found by computing
\begin{eqnarray}
\max \| Dx \|_2^2 \quad \st \| x \|_2 = 1.
\end{eqnarray}
ADMM is applied to the equivalent problem
\begin{eqnarray}
\min -\| Du \|_2^2 + \iota_{\{z:\, \| z\|_2 =1 \} }(v) \quad \st u-v=0,
\end{eqnarray}
where $\iota_{S}$ is the characteristic function defined by $\iota_{S}(v) = 0$, if $v\in S$, and $\iota_{S}(v) = \infty$, otherwise.

\section{Experiments \& Observations}
\begin{figure}[hbtp]
\centerline{
\includegraphics[width=0.5\linewidth]{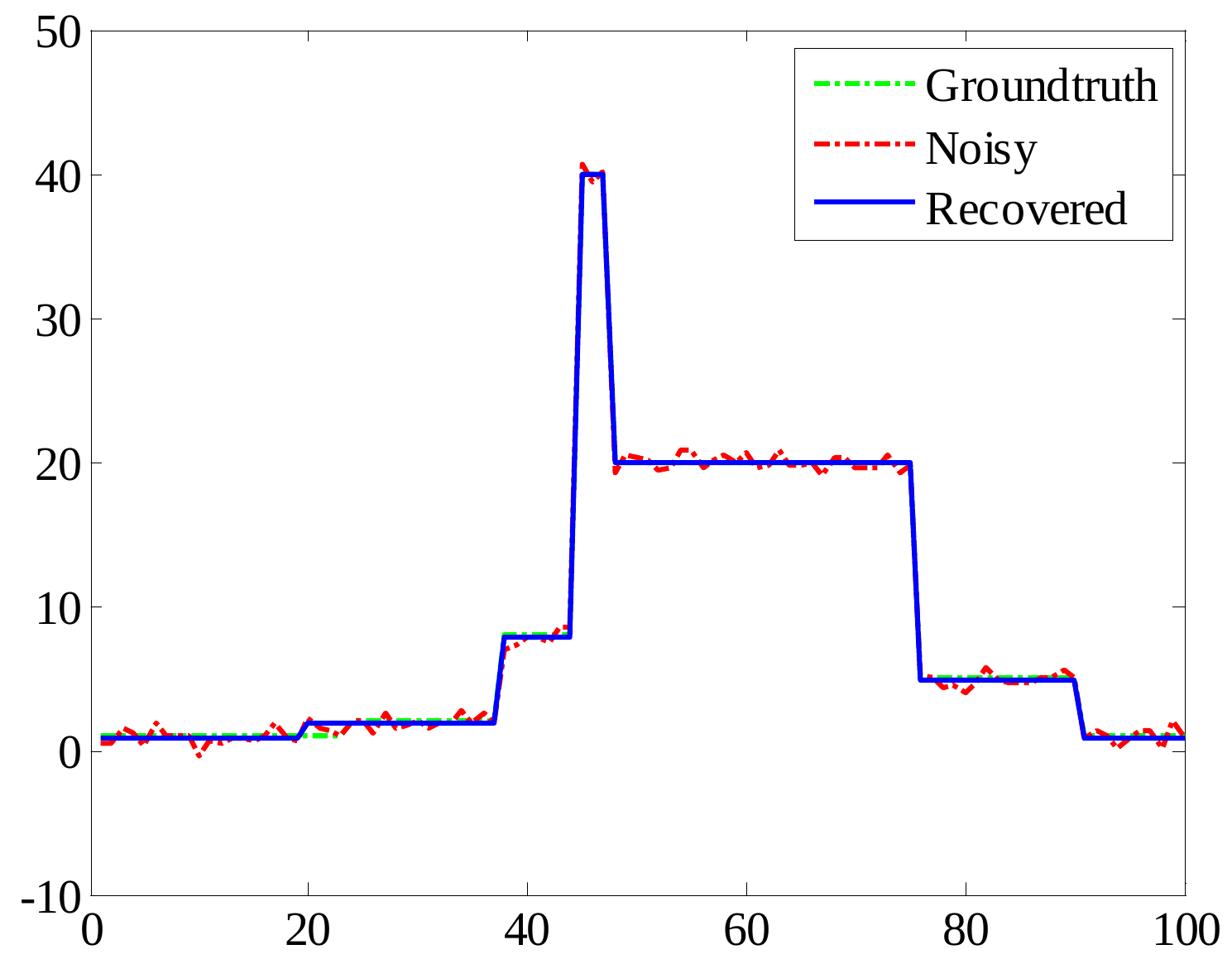}
}
\caption{\small The synthetic one-dimensional signal for $\ell_0$ regularized image denoising. The groundtruth signal, noisy signal (PSNR = 37.8) and recovered signal by AADMM (PSNR = 45.4) are shown. }
\label{fig:1dsig}
\end{figure}

\begin{figure}[hbtp]
\centerline{
\includegraphics[width=0.8\linewidth]{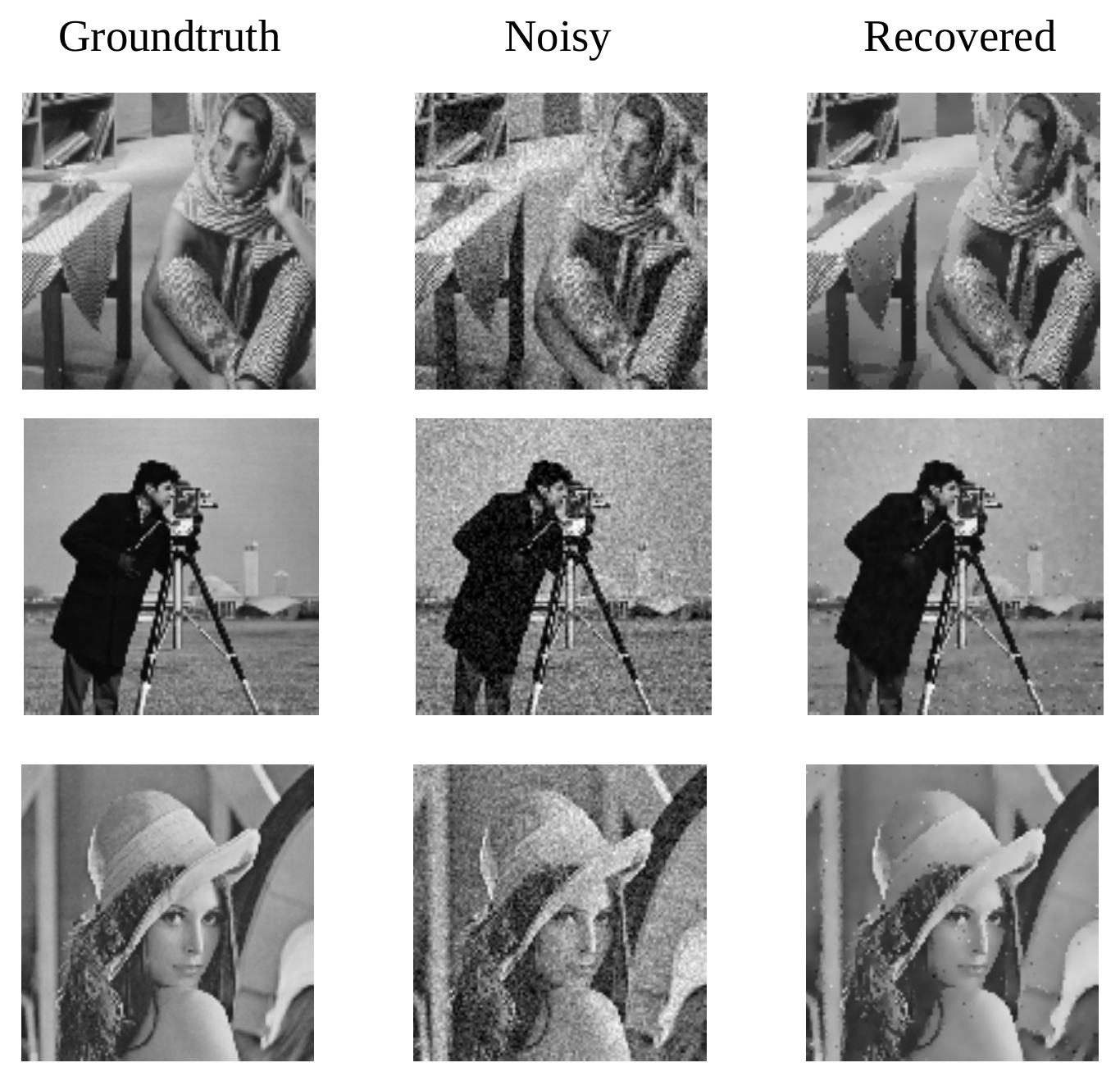}
}
\caption{\small The groundtruth image (left), noisy image (middle), and recovered image by AADMM (right) for $\ell_0$ regularized image denoising. The PSNR of the noisy/recovered images are 21.9/24.7 for "Barbara", 22.4/27.8 for "Cameraman", 21.9/27.9 for "Lena". }
\label{fig:l0tv_img}
\end{figure}

\textbf{Experimental setting.}
We implemented ``vanilla ADMM'' (ADMM with constant penalty),
and fast ADMM with Nesterov acceleration and restart \cite{goldstein2014fast}.  We also implemented two methods for automatically selecting penalty parameters:  residual balancing \cite{he2000alternating}, and the spectral adaptive method \cite{xu2016adaptive}. For $\ell_0$ regularized linear regression, the synthetic problem in \cite{zou2005regularization,goldstein2014fast,xu2016adaptive} and realistic problems in \cite{efron2004least,zou2005regularization,xu2016adaptive} are investigated with $\rho=1.$
For $\ell_0$ regularized image denoising, a one-dimensional synthetic problem was created by the process described in \cite{zou2005regularization}, and is shown in \cref{fig:1dsig}. For the total-variation experiments, the "Barbara" , "Cameraman", and "Lena" images are investigated, where Gaussian noise with zero mean and standard deviation 20 was added to each image (\cref{fig:l0tv_img}). $\rho =1$ and $\rho=500$ are used for the synthetic problem and image problems, respectively. For phase retrieval, a synthetic problem is constructed with a random matrix $D \in \bbC^{15000 \times 500}$, $x \in \bbC^{500}$, $e \in \bbC^{15000}$  and $c = \text{abs}(Dx + e)$. Three images in \cref{fig:l0tv_img} are used. Each image is measured with 21 octanary pattern filters as described in~\cite{candes2015phase}.
For the eigenvector problem, a random matrix $D \in \bbR^{20 \times 20}$ is used.

\begin{figure}[hbtp]
\centerline{
\includegraphics[width=\linewidth]{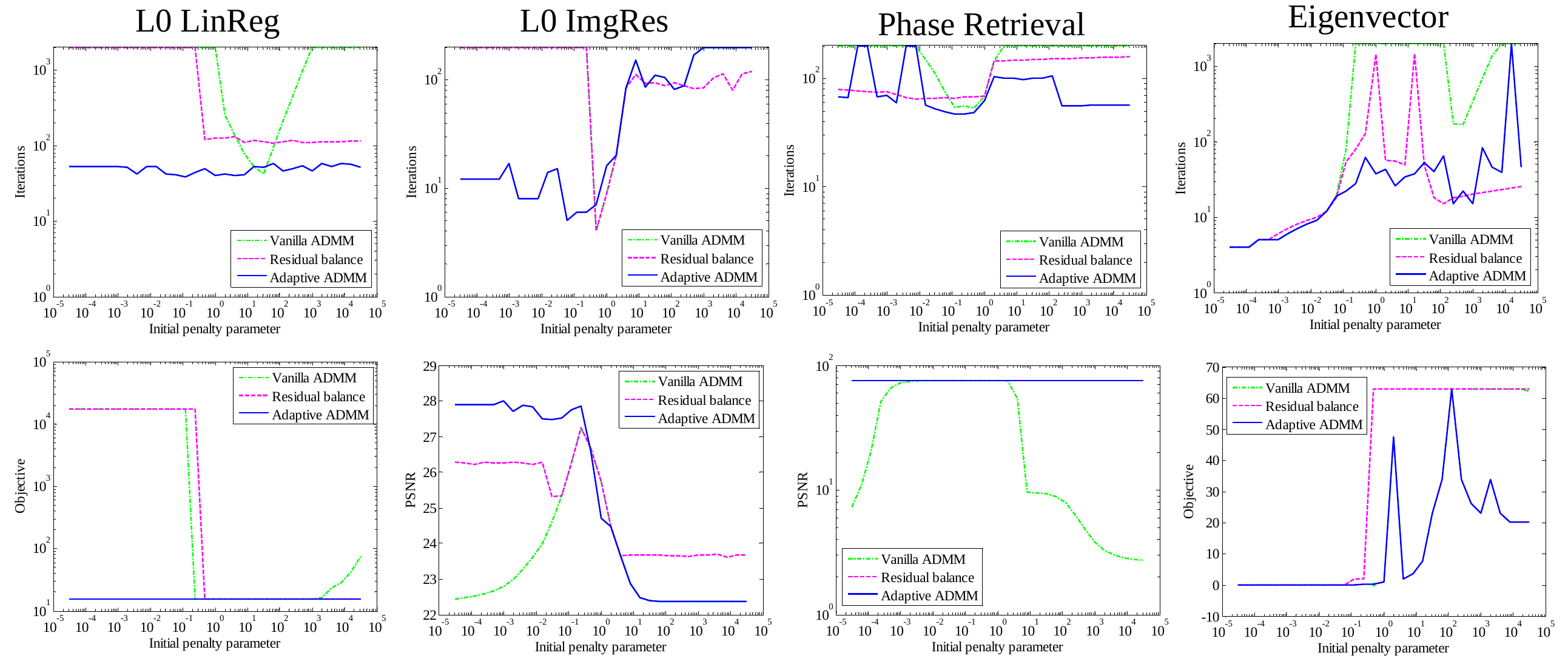}
}
\caption{\small Sensitivity to the (initial) penalty parameter $\tau_0$ for the $\ell_0$ regularized linear regression, eigenvector computation,  "cameraman" denoising, and phase retrieval. (top) Number of iterations needed as a function of initial penalty parameter. (bottom) The objective/PSNR of the minima found for each non-convex problem.
 }
\label{fig:tau}
\end{figure}

\textbf{Does ADMM converge in practice?}  The convergence of vanilla ADMM is quite sensitive to the choice of penalty parameter.  For vanilla ADMM, the iterates may oscillate, and if convergence occurs it may be very slow when the penalty parameter is not properly tuned.  The residual balancing method converges more often than vanilla ADMM, and the spectral adaptive ADMM converges the most often.  However, none of these methods uniformly beats all others, and it appears that vanilla ADMM with a highly tuned stepsize can sometimes outperform adaptive variants.


\textbf{Does the update order of $H(u)$ and $G(v)$ matter?} In \cref{fig:tau}, ADMM is performed by first minimizing with respect to the smooth objective term, and then the nonsmooth term. We repeat the experiments with the update order swapped, and report the results in \cref{fig:tau2} of the appendix. When updating the non-smooth term first, the convergence of ADMM for the phase retrieval problem becomes less reliable.  However, for some problems (like image denoising), convergence happened a bit faster than with the original update order.   Although the behavior of ADMM changes, there is no predictable difference between the two update orderings.

\newpage
\textbf{Is the local optimal solution good?} 
The bottom row of \cref{fig:tau} presents the objective/PSNR achieved by the ADMM variants when varying the (initial) penalty parameter.   In general, the quality of the solution depends strongly on the penalty parameter chosen.  There does not appear to be a predictable relationship between the best penalty for convergence speed and the best penalty for solution quality.

\textbf{Does the adaptive penalty work?}  In Table \ref{tab:exp}, we see that adaptivity not only speeds up convergence, but for most problem instances it also results in better minimizers.  This behavior is not uniform across all experiments though,  and for some problems a slightly lower objective value can be achieved using a finely tuned constant stepsize.

\begin{table}[htbp]
\centering
\caption{\small Iterations (with runtime in seconds) and objective (or PSNR) for the various algorithms and applications described in the text. Absence of convergence after $n$ iterations is indicated as $n+$.
}
\setlength{\tabcolsep}{3pt}
\begin{threeparttable}
\begin{tabular}{c|c|c||c|c|>{\bfseries}c}
\hline
Application & Dataset & \tabincell{c}{\#samples $\times$ \\ \#features\tnote{1}} & \tabincell{c}{Vanilla\\ ADMM}  & \tabincell{c}{Residual\\ balance~\cite{he2000alternating}} & \tabincell{c}{Adaptive \\ADMM \cite{xu2016adaptive}}\\
\hline\hline
\multirow{12}{*}{\tabincell{c}{$\ell_0$ regularized\\ linear regression}}
& \multirow{2}{*}{\tabincell{c}{Synthetic}} & \multirow{2}{*}{\tabincell{c}{50 $\times$ 40}} & 2000+(.621)  & 2000+(.604) & 39(.018) \\
&& & 1.71e4 & 1.71e4 & 15.2 \\
\cline{2-6}
& \multirow{2}{*}{\tabincell{c}{Boston}} & \multirow{2}{*}{\tabincell{c}{506 $\times$ 13}} & 2000+(.598) & 2000+(.570)  & 1039(.342) \\
&&& 1.50e5  & 1.50e5 & 1.34e5 \\
\cline{2-6}
& \multirow{2}{*}{\tabincell{c}{Diabetes}} & \multirow{2}{*}{\tabincell{c}{768 $\times$ 8}} & 2000+(.751)  & 2000+(.708)  & 28(.014) \\
&&& 384  & 648  & 285 \\
\cline{2-6}
&\multirow{2}{*}{\tabincell{c}{Leukemia}} & \multirow{2}{*}{\tabincell{c}{38 $\times$ 7129}} & 2000+(15.3)  & 78(.578)  & 63(.477) \\
&&& 19.0  & 19.0  & 19.0 \\
\cline{2-6}
& \multirow{2}{*}{\tabincell{c}{Prostate}} & \multirow{2}{*}{\tabincell{c}{97 $\times$ 8}} & 2000+(.413)  & 2000+(.466) & 29(.013) \\
&&& 1.14e3  & 380  & 324 \\
\cline{2-6}
&\multirow{2}{*}{\tabincell{c}{Servo}} & \multirow{2}{*}{\tabincell{c}{130 $\times$ 4}} &2000+(.426)  & 2000+(.471)  & 45(.014) \\
&&& 267  & 267 & 198 \\
\hline
\multirow{8}{*}{\tabincell{c}{$\ell_0$ regularized\\ image restoration}}
&\multirow{2}{*}{\tabincell{c}{Synthetic1D}} & \multirow{2}{*}{\tabincell{c}{100 $\times$ 1}} & 2000+(.701) & 1171(.409)  & 866(.319) \\
&&& 40.6 & 45.4 & 45.4 \\
\cline{2-6}
&\multirow{2}{*}{\tabincell{c}{Barbara}} & \multirow{2}{*}{\tabincell{c}{512 $\times$ 512}} & 200+(35.5)  & 200+(35.1)  & 18(3.33) \\
&&& 24.7  & 24.7  & 24.7 \\
\cline{2-6}
&\multirow{2}{*}{\tabincell{c}{Cameraman}}& \multirow{2}{*}{\tabincell{c}{256 $\times$ 256}} & 200+(5.75)  & 200+(5.60)  & 6(.190) \\
&&& 25.9  & 25.9  & 27.8 \\
\cline{2-6}
&\multirow{2}{*}{\tabincell{c}{Lena}}& \multirow{2}{*}{\tabincell{c}{512 $\times$ 512}} & 200+(35.5)  & 200+(35.8)  & 11(1.98) \\
&&& 25.9  & 25.9 & 27.9 \\
\hline
\multirow{8}{*}{\tabincell{c}{phase retrieval}}
&Synthetic & $15000 \times 500$ & 200+(19.4)  & 94(9.01) & 46(4.45) \\
\cline{2-6}
&\multirow{2}{*}{\tabincell{c}{Barbara}} & \multirow{2}{*}{\tabincell{c}{512 $\times$ 512 $\times$ 21}}  & 59(91.1) & 59(89.6) & 50(88.1) \\
&&& 81.5 & 81.5 & 81.5 \\
\cline{2-6}
&\multirow{2}{*}{\tabincell{c}{Cameraman}}& \multirow{2}{*}{\tabincell{c}{256 $\times$ 256 $\times$ 21}} & 59(29.6)  & 55(19.4)  & 48(20.8) \\
&&& 75.7  & 75.7 & 75.7\\
\cline{2-6}
&\multirow{2}{*}{\tabincell{c}{Lena}}& \multirow{2}{*}{\tabincell{c}{512 $\times$ 512 $\times$ 21}} & 59(90.1)) & 57(87.4) & 52(92.0) \\
&&& 81.4  & 81.5  & 81.5 \\
\hline
\end{tabular}%
\label{tab:exp}%
\begin{tablenotes}
    \item[1] width $\times$ height for image restoration; width $\times$ height $\times$ filters for phase retrieval
\end{tablenotes}
\end{threeparttable}
\end{table}%

\section{Conclusion}
We provide a detailed discussion of the performance of ADMM on several nonconvex applications, including $\ell_0$ regularized linear regression, $\ell_0$ regularized image denoising, phase retrieval, and eigenvector computation. In practice, ADMM usually converges for those applications, and the penalty parameter choice has a significant effect on both convergence speed and solution quality. Adaptive penalty methods such as AADMM \cite{xu2016adaptive} automatically select the penalty parameter, and perform optimization with little user oversight. For most problems, adaptive stepsize methods result in faster convergence or better minimizers than vanilla ADMM with a constant non-tuned penalty parameter.  However, for some difficult non-convex problems, the best results can still be obtained by fine-tuning the penalty parameter.


\clearpage

\small
\bibliographystyle{abbrvnat}
\bibliography{admm}
\clearpage

\section{Appendix: more experimental results}

\begin{figure}[hbtp]
\centerline{
\includegraphics[width=\linewidth]{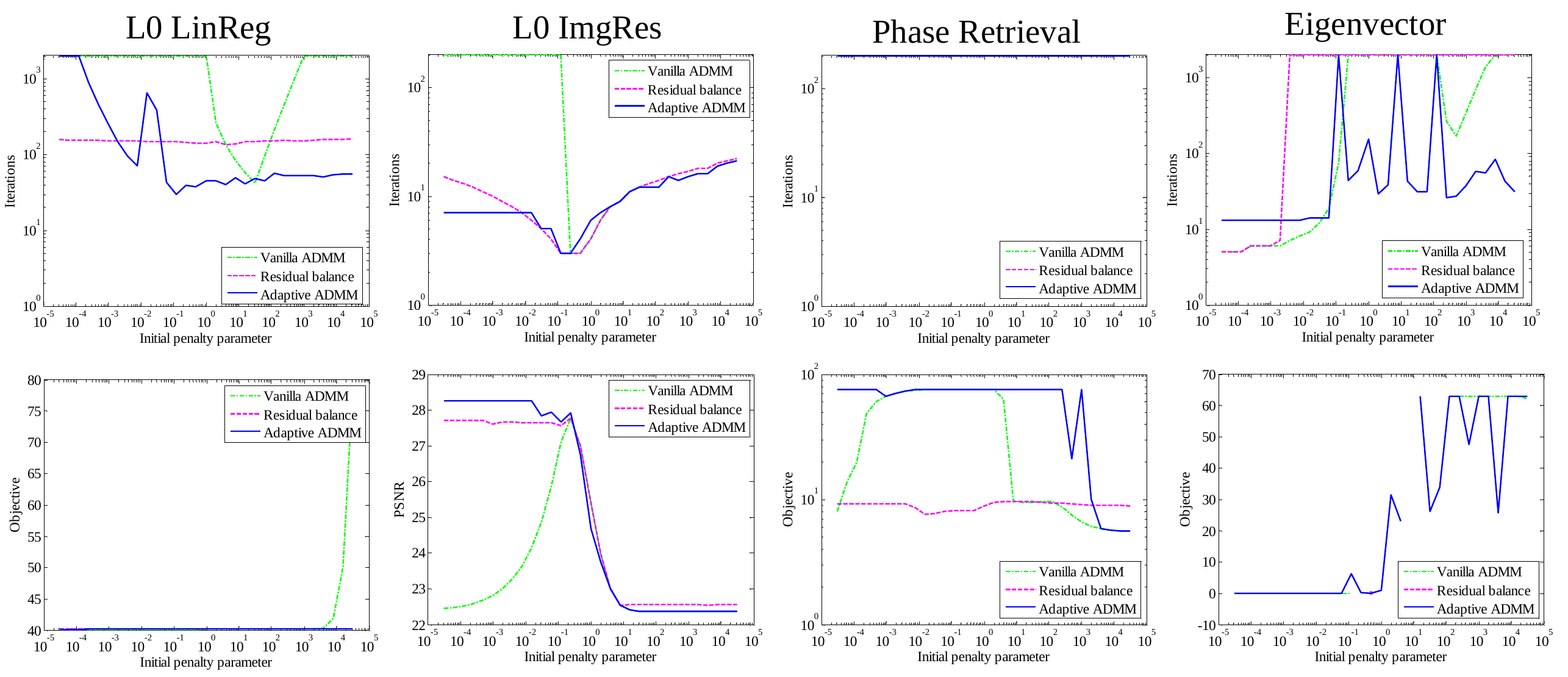}
}
\caption{\small Convergence results when the non-smooth objective term is updated first, and the smooth term is updated second. Sensitivity to the (initial) penalty parameter $\tau_0$ is shown for the synthetic problem of $\ell_0$ regularized linear regression, eigenvector computation, the  "cameraman"  denoising problem, and phase retrieval. The top row shows the convergence speed in iterations. The bottom row shows the objective/PSNR achieved by the final iterates.}
\label{fig:tau2}
\end{figure}

\section{Appendix: implementation details}

\subsection{$\ell_0$ regularized linear regression}
$\ell_0$ regularized linear regression is a nonconvex problem
\begin{eqnarray}
\min_x \frac{1}{2} \| Dx - c \|_2^2 + \rho \|x\|_0 \label{eq:l0r}
\end{eqnarray}
where $D \in \bbR^{n \times m}$ is the data matrix, $c$ is the measurement vector, and $x$ is the regression coefficients. ADMM is applied to solve problem~\eqref{eq:l0r} by solving the equivalent problem
\begin{eqnarray}
\min_{u,v} \frac{1}{2} \| Du - c \|_2^2 + \rho \|v\|_0~~~~\st~~u-v=0.
\end{eqnarray}
The proximal operator of the $\ell_0$ norm is the hard-thresholding,
\begin{eqnarray}
\eqhd(z, t) = \arg\min_x \|x\|_0 + \frac{1}{2t} \| x - z \|_2^2  =  z \odot \mcI_{\{z:|z| > \sqrt{2t}\}}(z),
\end{eqnarray}
where $\odot$ represents element-wise multiplication,  and $\mcI_{S}$ is the indicator function of the set $S$: $\mcI_{S}(v) = 1$, if $v\in S$, and $\mcI_{S}(v) = 0$, otherwise.
Then the steps of ADMM can be written
\begin{eqnarray}
u\kp &=& \arg\min_{u} \| Du - c \|_2^2 + \frac{\tau}{2} \| 0-u+v\itk + \gl\itk/\tau\|_2^2 \\
&=&
\begin{cases}
(D^{T}D+\tau I_n)^{-1}(\tau v\itk+\gl\itk + D^{T}c) &\eqif~~ n \geq m\\
(I_n - D^{T} (\tau I_m + DD^{T})^{-1} D)(v\itk+\gl\itk/\tau + D^{T}c/\tau) & \eqif~~ n<m
\end{cases}\\
v_{k+1} &=& \arg\min_{v} \rho \|v\|_0 + \frac{\tau}{2} \| 0 -u\kp +v + \gl\itk/\tau\|_2^2 = \eqhd( u\kp-\gl\itk/\tau, \rho/\tau) \\
\gl\kp &=& \gl\itk +\tau (0 -  u\kp + v\kp).
\end{eqnarray}

\subsection{$\ell_0$ regularized image denoising}
The $\ell_0$ regularizer~\cite{dong2013efficient} is an alternative to the $\ell_1$ regularizer when computing total variation~\cite{goldstein2009split,goldstein2014fast}. $\ell_0$ regularized image denoising solves the nonconvex problem
\begin{eqnarray}
\min_{x} \frac{1}{2} \| x-c \|_2^2 + \rho \| \grad x\|_0
\end{eqnarray}
where $c$ represents a given noisy image, $\grad$ is the linear gradient operator, and $\| \cdot \|_2$/$\| \cdot \|_0$ denotes the $\ell_2$/$\ell_0$ norm of vectors. The steps of ADMM for this problem are
\begin{align}
u\kp &= \arg\min_u \frac{1}{2} \| u-c \|_2^2 + \frac{\tau}{2} \|  v\itk + \gl\itk/\tau - \grad u\|_2^2\\ & = (I + \tau \grad^T\grad)^{-1} (c+\tau\grad^T(v\itk+\gl\itk/\tau))\\
v\kp &= \arg\min_{v} \rho \|v\|_0 + \frac{\tau}{2} \| 0 - \grad u\kp +v + \gl\itk/\tau\|^2 = \eqhd( \grad u\kp -\gl\itk/\tau, \rho/\tau) \\
\gl\kp &= \gl\itk +\tau (0 -  \grad u\kp + v\kp)
\end{align}
where the linear systems can be solved using fast Fourier transforms.

\subsection{Phase retrieval}
Ptychographic phase retrieval~\cite{yang2011iterative,wen2012alternating} solves problem
\begin{eqnarray}
\min_{x} \frac{1}{2}||\text{abs}(Dx)-c||_2^2,
\end{eqnarray}
where $x\in \bbC^n$, $D\in \bbC^{m \times n}$, and $\text{abs}(\cdot)$ denotes the elementwise magnitude of a complex-valued vector. ADMM is applied to the equivalent problem
\begin{eqnarray}
\min_{u,v} \frac{1}{2}||\text{abs}(u)-c||_2^2 ~~\st u-Dv = 0.
\end{eqnarray}
Define the projection operator of a complex valued vector as
\begin{eqnarray}
\text{absProj}(z, c, t) = \min_x \frac{1}{2} \| \text{abs}(x)-c \|_2^2 + \frac{t}{2} \| x-z\|_2^2 = \left(\frac{t}{1+t}\,\text{abs}(z) + \frac{1}{1+t}\,c\right) \odot \text{sign}(z),
\end{eqnarray}
where $\text{sign}(\cdot)$ denotes the elementwise phase of a complex-valued vector.
In the following ADMM steps, notice that the dual variable $\gl \in \bbC^m$ is complex, and the penalty parameter $\tau \in \bbR$ is a real non-negative scalar,
\begin{align}
u\kp &= \arg\min_u \frac{1}{2} \| \text{abs}(u)-c \|_2^2 + \frac{\tau}{2} \|   Dv\itk + \gl\itk/\tau - u\|_2^2
=  \text{absProj}(Dv\itk + \gl\itk/\tau, c, \tau)\\
v\kp &= \arg\min_v \, 0 + \frac{\tau}{2} \| 0 -  u\kp + Dv + \gl\itk/\tau\|_2^2 = D^{-1}( u\kp - \gl\itk /\tau) \\
\gl\kp &= \gl\itk +\tau (0 -   u\kp + D v\kp).
\end{align}

\subsection{Eigenvector problem}
The eigenvector problem is a fundamental problem in numerical linear algebra.  The leading eigenvector of a matrix can be recovered by solving the Rayleigh quotient maximization problem
\begin{eqnarray}
\max \| Dx \|_2^2 \quad \st \| x \|_2 = 1.
\end{eqnarray}
ADMM is applied to the equivalent problem
\begin{eqnarray}
\min -\| Du \|_2^2 + \iota_{\{z:\, \| z\|_2 =1 \} }(v) \quad \st u-v=0,
\end{eqnarray}
where $\iota_{S}$ is the characteristic function of the set $S$: $\iota_{S}(v) = 0$, if $v\in S$, and $\iota_{S}(v) = \infty$, otherwise. The ADMM steps are
\begin{align}
u\kp &= \arg\min_u -\| Du \|_2^2 + \frac{\tau}{2} \| 0 -  u +v\itk + \gl\itk/\tau\|_2^2
= (\tau I - 2 D^TD)^{-1} (\tau v\itk + \gl\itk)\\
v\kp &= \arg\min_{v} \iota_{\{z:\, \| z\|_2 =1 \} }(v) + \frac{\tau}{2} \| 0 -  u\kp +v + \gl\itk/\tau\|^2 = \frac{u\kp - \gl\itk/\tau}{\| u\kp - \gl\itk/\tau \|_2}\\
\gl\kp &= \gl\itk +\tau (0 -  u\kp + v\kp).
\end{align}

\section{Appendix: synthetic and realistic datasets}

We provide the detailed construction of the synthetic dataset for our linear regression experiments. The same synthetic dataset has been used in~\cite{zou2005regularization,goldstein2014fast,xu2016adaptive}. Based on three random normal vectors $\nu_a, \nu_b, \nu_c \in \bbR^{50}$, the data matrix $D = [d_1 \ldots d_{40}] \in \bbR^{50 \times 40}$ is defined as
\begin{eqnarray}
d_i =
\begin{cases}
\nu_a + e_i, &i=1,\ldots,5, \\
\nu_b + e_i, &i=6,\ldots,10, \\
\nu_c + e_i, &i=11, \dots,15, \\
\nu_i\in N(0,1), &  i=16, \dots,40,
\end{cases}
\end{eqnarray}
where $e_i$ are random normal vectors from $N(0, 1)$. The problem is to recover the vector
\begin{eqnarray}
x^* =
\begin{cases}
3, &i=1,\ldots,15,\\
0, &\text{otherwise}
\end{cases}
\end{eqnarray}
from noisy measurements of the form
$
c = Dx^* + \hat{e}, $
 with $\hat{e} \in N(0, 0.1)$

\end{document}